\documentclass[11pt,twoside,a4paper]{article}
\usepackage{amsmath,amssymb,amsthm}

\newtheorem{Thm}{Theorem}[section]
\newtheorem{Lem}[Thm]{Lemma}
\newtheorem{Pro}[Thm]{Proposition}
\newtheorem{Cor}[Thm]{Corollary}
\newtheorem{Def}[Thm]{Definition}

\theoremstyle{definition}

\theoremstyle{remark}
\newtheorem{Rem}[Thm]{Remark}

\newcommand{\R}{\mathbb{R}}
\newcommand{\Z}{\mathbb{Z}}
\newcommand{\N}{\mathbb{N}}
\renewcommand{\H}{\mathbb{H}}

\renewcommand{\phi}{\varphi}

\newcommand{\rank}{\operatorname{rank}}

\newcommand{\cone}{\operatorname{Cone}}

\newcommand{\lev}{\operatorname{lev}}
\newcommand{\fin}{\operatorname{fin}}

\newcommand{\wh}{\widehat}

\begin{document}

\title{Embedding of Coxeter groups in a product of trees}
\author{Alexander Dranishnikov\footnote{Supported by
NSF} \ \& Viktor Schroeder\footnote{Partially supported by SNF}}

\maketitle

\begin{abstract}
We prove that a right angled Coxeter group $\Gamma$ with chromatic number
$n$ can be embedded in a bilipschitz way 
into the product of n locally finite trees.
We give applications of this result to various embedding problems
and determine the hyperbolic rank of products of exponentially branching trees.
\end{abstract}

\section{Introduction}

We consider a finitely generated right angled Coxeter group $\Gamma$ , 
i.e. a group
$\Gamma$ together with a finite set of 
generators $S$, such that every element of $S$ has order two and that all
relations in $\Gamma$ are consequences of relations of the form
$st = ts$, where $s,t \in S$.

We prove embedding results of the Cayley graph 
$C(\Gamma ,S)$ into products of
trees. On graphs and trees we consider always the simplicial metric, hence
every edge has length 1. On a product of trees we consider the
$l_1$-product metric, i.e. the distance is equal to the sum of
the distances in the factors.

In [DJ] it was shown that the Cayley graph of a Coxeter group
admits an equivariant isometric embedding  into a finite product of 
locally infinite trees.
Here we give a better estimate on the number of factors in the 
right-angled case.
The estimate is given in terms of the chromatic number.
Consider therefore  {\it colourings }$c: S \to \{1,\ldots ,n\}$ 
with the property that for different
$s,t \in S$ with $st = ts$ we have $c(s) \neq c(t)$. The minimal 
number $n$ of colours
needed is called the chromatic number of $\Gamma$.

\begin{Thm} \label{tree}
Suppose that the chromatic number of a right-angled 
Coxeter 
group $\Gamma$ is $n$. Then the Cayley graph $C(\Gamma,S)$ 
admits
an equivariant isometric embedding into the product of 
$n$
simplicial trees.
\end{Thm}

Besides of trivial cases, these trees are  locally infinite.
However we are able to embed the Cayley graph bilipschitz into
a product of locally compact trees.

\begin{Def}

A pointed simplicial tree $(T,t_0)$ is called exponentially 
branching, if there exists a
number $\sigma >0$ such that every vertex $t \in T$ has more
than $e^{\sigma d(t,t_0)}$ neighbours where $d$ is the metric on $T$.

\end{Def}

\begin{Thm} \label{locallycompacttree}
Let $\Gamma$ be a right-angled 
Coxeter group with chromatic number $n$, let $T$ be an exponentially branching
locally compact simplicial tree, and let $r > 0$ be a number.
 Then there exists  bilipschitz embedding 
$\psi : C(\Gamma,S) \to T\times \cdots \times T$ ($n$-factors),
such that $\psi$ restricted to every ball of radius $ r$ is isometric.

\end{Thm}

It is an interesting open problem, if a corresponding embedding result 
holds for trees with bounded valence.

We can apply Theorem \ref{locallycompacttree} for a special
Coxeter group operating on the hyperbolic plane $\H^2$ and obtain:

\begin{Cor} \label{hypplane}
For every exponentrially branching tree $T$ there exists a bilipschitz embedding
$\phi: \H^2 \rightarrow T\times T$.
\end{Cor}

Combining this with a result of Brady and Farb we get the following higher
dimensional version:

\begin{Cor} \label{higherdim}
For every exponentrially branching tree $T$ there exists a bilipschitz embedding
$\psi: \H^n \rightarrow T\times\cdots\times T$ of the 
hyperbolic space $\H^n$ into the $2(n-1)$ fold product of $T$.
\end{Cor}

It is an open question, if for 
$n \geq 3$ there is a bilipschitz embedding of
$\H^n$ into the $n$-fold product of locally compact trees.
There are two partial results in this direction. In
\cite{BS2} it is show that there exists a quasiisometric embedding of
$\H^n$ into an $n$-fold product of locally infinite trees.
On the other hand
a recent construction of 
Januszkiewicz and 
Swiatkowski \cite{JS} shows for every $n$ the existence of
a right angled Gromov hyperbolic Coxeter group with virtual
cohomological dimension and colouring number equal to $n$.
Combining Theorem \ref{locallycompacttree} with that result
we obtain:

\begin{Cor} \label{vcohdim}

For every exponentially branching tree $T$ and any given number
$n$ there exits a Gromov hyperbolic group $\Gamma_n$ with virtual
cohomological dimension $n$ and a bilipschitz embedding of
the Cayley graph of $\Gamma_n$ into the product
$T\times \cdots \times T$ ($n$-factors).

\end{Cor}

Corollary \ref{vcohdim} can be used to determine the hyperbolic rank
(compare \cite{BS1}) of a product of trees:

\begin{Cor} \label{hyprank}
The hyperbolic rank of the product 
of
$n$ trees with exponential branching is $(n-1)$.
\end{Cor}

%%%%%%%%%%%%%%%%%%%%%%%%%%%%%%%%%%%%%%%%%%%%%%%%%%%%%%%%%%%%%%%%%%%%%%%%%%%%%%%%%

\section{Right Angled Coxeter Groups}

In this section we review the necessary facts from the theory of
right angled Coxeter groups.

\subsection{Deletion Rule and Cayley Graph}

A Coxeter matrix $(m_{s,t})_{s,t \in S}$ is a symmetric
$S \times S$ matrix with 1 on the diagonal and with all other
entries nonnegative integers different from 1.
A Coxeter matrix defines a Coxeter group $\Gamma$ generated by
the index set $S$ with relations
$(st)^{m_{s,t}} = 1$ for all $s,t \in S$.
Here we use the convention that $\gamma^0 =1$ for all elements, 
thus if $m_{s,t}=0$ then there
is no relation between $s$ and $t$. A Coxeter group $\Gamma$ 
is {\em finitely generated}, if
$S$ is finite. The group $\Gamma$ is called {\em right angled}, if all
entries of the corresponding Coxeter matrix are 0,1,2. The 
Coxeter matrix of a right
angled Coxeter group is completely described by a graph 
with vertex set $S$ where we
connect two vertices $s$ and $t$ iff $m(s,t)=2$.

Right angled Coxeter groups have a very simple deletion law. By the following
two operations every word $w$ in the generators $S$ can be transformed 
to a reduced word and two reduced
words representing the same element can be transformed by means only the second
operation \cite{Br} :

\begin{description}
\item[(i)] delete a subword of the form $ss$, $s\in S$

\item[(ii)] replace a subword $st$ by $ts$ if $m_{s,t}=2$.
\end{description}

\noindent This deletion rule has the following consequences:

\begin{Lem} \label{deletionrule}

\begin{description}
\item[(a)] If $w$ and $w'$ represent the same element then the lenghts of
$w$ and $w'$ are either both even or both odd.

\item[(b)] Let $w$ and $w'$ be reduced representations of the same element
$\gamma \in \Gamma$, then $w$ and $w'$ are formed from the same set of
letters and they have the same length.

\end{description}

\end{Lem}

We now investigate some properties of the Cayley graph
$C(\Gamma,S)$ of $\Gamma$ with respect to the generating set $S$ for
a right angled Coxeter group.
Let $\gamma \in \Gamma$ and let $w$ be a reduced representation of
$\gamma$. The length of $w$ is denoted by
$\ell (\gamma )$ and called the norm of $\gamma$. This is well defined by (b).
On $\Gamma$ we consider the distance function
$d(\gamma,\beta) =\ell (\gamma^{-1} \beta )$.
Let $\gamma$ and $\beta$ be elements of $\Gamma$ which are
neighbours in the Cayley graph and let $w$ be a reduced word representing
$\gamma$. Then there exists a generator $s \in S$ such that
$\beta$ has the representation $ws$. It follows from (b) that
$\ell (\gamma ) \neq \ell (\beta )$. Thus an edge in the Cayley graph
connects two elements with different norm.
This allows us to orient the edges of the Cayley 
graph and we choose the orientation in the
way that the initial point of an edge is closer to the 1-element than
the endpoint of the edge.
We write an oriented edge as
$[\alpha,\beta]$.
As usual one can define geodesics in the Cayley graph.
A geodesic between two points $\alpha,\beta\in\Gamma$ is given by a
sequence
$\alpha=\gamma_0,\ldots,\gamma_k=\beta$
with $d(\gamma_i,\gamma_j)=\ \mid i-j\mid$.

Let $\alpha, \beta$ and $\gamma$ be elements of $\Gamma$. We say that
$\gamma$ {\it lies between} $\alpha$ {\it and} $\beta$ if
$d(\alpha,\gamma)  + d(\gamma,\beta) = d(\alpha,\beta)$.
The Cayley graphs of right angled Coxeter groups have the
following property, which says that any three
points in $\Gamma$ span a tripoid:

\begin{Lem} \label{between}
Let $\alpha,\beta,\gamma \in \Gamma$, then there exists $\delta \in \Gamma$
such that $\delta$ lies between $\gamma_i$ and $\gamma_j$ for
any choice of distinct elements $\gamma_i,\gamma_j \in \{ \alpha,\beta,\gamma\}$.
\end{Lem}

\begin{proof}

By the $\Gamma$-invariance of the metric it suffices to show this result
for the case that $\gamma = 1$.
Let $\alpha,\beta \in \Gamma$ and consider a geodesic path
$\alpha =\alpha_0,\ldots,\alpha_k=\beta$ from $\alpha$ to
$\beta$ and consider the sequence of norms
$n_0=\ell (\alpha_0 ),\ldots,n_k =\ell (\alpha_k )$. Note that 
by the properties discussed
above $\mid n_i - n_{i+1}\mid\ =1$.
Assume that there is a subsequence $\alpha_{i-1},\alpha_i,\alpha_{i+1}$ with
$n_{i-1}<n_i>n_{i+1}$.
Then one can represent $\alpha_i$ in a reduced way as $w_1t = w_2s$ where
$w_1$ is a reduced word representing $\alpha_{i-1}$ and
$w_2$ a reduced word representing $\alpha_{i+1}$ and $s,t \in S$.
Since by the deletion law one can transform $w_1t$ into $w_2s$ by means
of operations of type (ii), we see that $st=ts$ and that one can represent
$\alpha_i$ by a reduced word of the form $wst=wts$ where $ws$ represents
$\alpha_{i-1}$ and $wt$ represents $\alpha_{i+1}$. Replace now $\alpha_i$ by
the element $\alpha_i'$
represented by $w$ and we obtain a new geodesic sequence between $\alpha$ and
$\beta$ such that for the corresponding sequence of norms
we have
$n_{i-1}>n_i'<n_{i+1}$.
Applying this procedure several times
we obtain a geodesic path from $\alpha$ to $\beta$, such that the sequence of the norms
has no local maximum any more and hence only a global minimum
$n_{i_0}$. The corresponding element $\delta = \alpha_{i_0}$ lies
between $\alpha$ and $\beta$ but also between
$1$ and $\alpha$ resp. $1$ and $\beta$.

\end{proof}

Let $s \in S$ be a generator.
We define $Z_s(S)=\{t \in S \mid st=ts\}$ and 
$Z_s^0(S)= Z_s(S)\setminus \{s\}$.
By $Z_s(\Gamma)$ we denote the centralizer of $s$ (in $\Gamma$).

\begin{Lem} 
For a 
right-angled Coxeter group the centralizer $Z_s(\Gamma)$ of a 
generator
$s\in S$ is the subgroup generated by the set $Z_s(S)\subset S$.
\end{Lem}

\begin{proof}
Clearly all elements in $Z_s(S)$ commute with $s$. Suppose that
$w$ commutes with $s$. Then $wsw^{-1}=s$. Let $w=s_1\dots s_k$ be a
reduced presentation. Then $s_1\dots s_kss_k\dots s_1=s$. Note 
that $s_k$ commutes
with $s$. If it does not commute, the letter $s$ in the middle
cannot be canceled and hence after every transformation 
there always will be a letter $s_k$ to the left from $s$. This 
contradicts that the
length of the word is one. If $s_k$ commutes with $s$ it can be canceled and we
consider $s_{k-1}$. By induction we have that all $s_i$ commute with $s$.
Therefore $s_i\in Z_s(S)$ for $i=1,\dots,k$.

\end{proof}

If $w$ is a reduced representation of an element in $Z_s(\Gamma)$ 
then all letters
of $w$ commute with $s$ and hence the letter $s$ can only occur once.
It follows that $Z_s(\Gamma)$ spits in a natural way as
$Z_s^0(\Gamma)\times \Z_2$, where
$Z_s^0(\Gamma)$ is the subgroup of $\Gamma$ generated by $Z_s^0(S)$.
Note that $Z_s^0(\Gamma)$ is itself a right angled Coxeter group
with generating set $Z_s^0(S)$.

\begin{Lem} \label{reduced}

Let $s_1\ldots s_m$ be a reduced word and $s \in S$.
Then the following equivalence holds:
$s s_1\ldots s_m$ is not reduced $\Leftrightarrow
\exists \ i\in\{ 1,\ldots,m\}$ such that $s_i = s$ and
$s_1,\ldots,s_{i-1} \in Z^0_s(S)$.

\end{Lem}

\begin{proof}

$" \Leftarrow "$  is clear.

$" \Rightarrow "$
By assumption the word $s s_1\ldots s_m$ is not reduced, let
$s_1',\ldots,s_l'$ be a reduced representation of
$s s_1\ldots s_m$. It follows from Lemma \ref{deletionrule} that
$l=m-r$ where $r$ is a positive odd number.
Since $s s_1'\ldots s_l'$ is a presentation of
$s_1\ldots s_m$ we have $l+1 \geq m$ and thus
$l=m-1$. Hence $s s_1'\cdots s_l'$ is a minimal
representation of $s_1\ldots s_m$. By Lemma \ref{deletionrule} (b)
the letter $s$ occurs in
$\{ s_1\ldots s_m\}$.
Let $i$ be the smallest integer such that $s_i=s$.
Assume that $i \geq 2$ and $s_1$ not in $Z_s^0(S)$.
Then the letter $s$ in front of $s_1\ldots s_m$
can not be cancelled with another letter $s$ in
$s_1\ldots s_m$. Since $s_1\ldots s_m$ is reduced
this implies that also $s s_1\ldots s_m$ is reduced.
This contradiction shows that
$s_1 \in Z_s^0(S)$. Inductively it follows that
$s_1,\ldots,s_{i-1} \in Z^0_s(S)$.

\end{proof}

\begin{Lem}

Let $s,t \in S$ and let $s_1\ldots s_m$ be a reduced word,
such that the words
$s_1\ldots s_m t$ and $s s_1\ldots s_m$ are reduced, but
the word $s s_1\ldots s_m t$ is not reduced.
Then $s=t$ and $s_1,\ldots,s_m \in Z^0_s(S)$.
\end{Lem}

\begin{proof}

Define $s_{m+1} := t$.
Since $s_1\ldots s_m s_{m+1}$ is reduced and
$s s_1\ldots s_m s_{m+1}$ is not reduced, it follows from
Lemma \ref{reduced}
that there exists
$i\in\{1,\ldots,m+1\}$ such that
$s_i =s$ and
$s_1,\ldots, s_i \in Z^0_s(S)$.
Since $s_1\ldots s_m$ is reduced and $s s_1\ldots s_m$ is
also reduced, it follows from Lemma \ref{reduced} that
$i= m+1$ which implies that
$s=t$ and $s_1,\ldots,s_m \in Z^0_s(S)$.

\end{proof}

\noindent For a generator $s \in S$ we define
the {\it halfspace}
$H_s=\{\gamma \in \Gamma \mid \ell (s\gamma ) > \ell (\gamma) \}$
with the boundary
$\partial H_s =\{\gamma\in H_s \mid d(\gamma,\Gamma \setminus H_s) =1\}$ .
We have the following properties

\begin{Lem} \label{halfspace}

\begin{description}

\item[(a)] $\Gamma = H_s \cup sH_s$ is a disjoint union.

\item[(b)] If $\alpha \in H_s$ and $\alpha t \in sH_s$ for some $t\in S$.
Then $\ell (\alpha )<\ell (\alpha t )$, $t=s$ and $\alpha \in Z^0_s(\Gamma)$.

\item[(c)] $\partial H_s =  Z_s^0(\Gamma)$ and $H_s$ is 
$Z_s^0(\Gamma)$-invariant.

\item[(d)] $H_s$ and $\partial H_s$ are totally convex, i.e. every 
geodesic with
initial and endpoint in $H_s$ (resp. in $\partial H_s$ ) is 
completely contained in $H_s$
(resp. in $\partial H_s$).

\end{description}
\end{Lem}

\begin{proof}

(a)  Since trivially $\ell (s\gamma ) > \ell (\gamma )$ if and only if
$\ell (ss\gamma ) < \ell (s\gamma )$ we have $H_s \cap sH_s = \emptyset$.
For given $\gamma \in \Gamma$ we know that
$\ell (\gamma ) \neq \ell (s \gamma )$. 
Thus $\gamma \in H_s$ or $s\gamma \in H_s$ which
implies $\Gamma = H_s \cup sH_s$.

(b) Assume to the contrary that $\ell (\alpha ) >\ell (\alpha t )$.
Let $s_1\ldots s_k$ be a reduced representation of $\alpha t$,
then $s_1\ldots s_k t$ is a reduced representation of
$\alpha $. Since $s s_1\ldots s_k$ is not reduced it follows
trivially that $s s_1\ldots s_k t$ is also not reduced.
This is a contradiction to $\alpha \in H_s$.
Thus  $\ell (\alpha ) <\ell (\alpha t )$.

Let now $s_1\ldots s_k$ be a reduced representation of $\alpha$
then $s_1\ldots s_k t$ is a reduced representation of
$\alpha t$. By assumption $s s_1\ldots s_k$ is reduced and
$s s_1\ldots s_k t$ is not reduced.
Then $t=s$ and $\alpha \in Z^0_s(\Gamma)$ by Lemma \ref{reduced}.

(c) $\partial H_s \subset Z_s^0(\Gamma)$ follows immediately from (b).
If $s_1\ldots s_k$ is a reduced representation of an element 
$\gamma \in Z^0_s(\Gamma)$ with $s_i \in Z_s^0(S)$, then
$s s_1\ldots s_k$ is a reduced word, i.e. $Z_s^0(\Gamma) \subset H_s$.
Clearly $s s_1\ldots s_k \in sH_s$ thus $Z_s^0(\Gamma) \subset \partial H_s$.

To show that $H_s$ is  $Z_s^0(\Gamma)$-invariant, we first prove
that $H_s$ is star shaped with respect to $1$,
i.e. every geodesic from $1$ to $\gamma \in H_s$ is completely
contained in $H_s$. Such a geodesic corresponds to a reduced
representation  $s_1\ldots s_m$ of $\gamma$.
Since $\gamma \in H_s$ we see that
$s s_1\ldots s_m $ is a reduced word and hence
$s s_1\ldots s_k $ is reduced for all $0\leq k\leq m$ which
implies that the geodesic lies in $H_s$.

Let now $\gamma \in H_s$ and $\alpha \in Z_s^0(\Gamma)$.
Let $1=\beta_0,\ldots,\beta_m=\gamma$ be a geodesic
which is by the above completely contained in $H_s$.
Then $\alpha \beta_0,\ldots,\alpha\beta_m$ is a geodesic
from $\alpha \in H_s$ to $\alpha\gamma$.
If $\alpha\gamma$ is not contained in $H_s$, then
this geodesic leaves $H_s$ and there is $i$ such
that $\alpha\beta_i \in H_s$ and $\alpha \beta_{i+1} \in sH_s$.
By (b) $\alpha\beta_{i} \in Z_s^0(\Gamma)$ and
$\alpha\beta_{i+1} \in Z_s(\Gamma)\setminus Z_s^0(\Gamma)$.
Since $\alpha\in Z_s^0(\Gamma)$ this implies
$\beta_{i+1} \in (Z_s(\Gamma)\setminus Z_s^0(\Gamma)) \subset sH_s$,
a contradiction.

(d) Let $\alpha_0,\ldots,\alpha_m$ be a geodesic with initial
and endpoint in $H_s$. If this geodesic is not completely
contained in $H_s$ let $i$  be the smallest index
such that $\alpha_i \in sH_s$ and $j$ be the largest index
with $\alpha_j \in sH_s$.
By (b) $\alpha_{i-1},\alpha_{j+1} \in Z_s^0(\Gamma)$.
Thus it remains to show that every geodesic
joining two points $\alpha,\beta\in Z_s^0(\Gamma)$ is
completely contained in $H_s$.
A geodesic  from $1$ to $\alpha^{-1}\beta \in Z_s^0(\Gamma)$
corresponds to a reduced word representing $\alpha^{-1}\beta$.
By the deletion rule such a word is formed only out of letters
from $Z_s^0(S)$.
Thus every geodesic from
$1$ to $\alpha^{-1}\beta$ is contained in $Z_s^0(\Gamma)\subset H_s$ and
since $H_s$ is $\alpha$-invariant by (c) every geodesic from
$\alpha$ to $\beta$ is contained in $H_s$.
The argument also shows that every geodesic between two
points $\alpha ,\beta \in Z_s^0(\Gamma)=\partial H_s$
is conpletely contained in $\partial H_s$.

\end{proof}

\begin{Rem}
Lemma \ref{halfspace} says that the boundary of 
$\partial H_s =\{\gamma\in H_s, d(\gamma,\Gamma \setminus H_s) =1$ 
is equal to 
$Z_s^0(\Gamma)$ and $\partial sH_s = s Z_s^0(\Gamma)$.
\end{Rem}

%%%%%%%%%%%%%%%%%%%%%%%%%%%%%%%%%%%%%%%%%%%%%%%%%%%%%%%%%%%%%%%
%%%%%%%%%%%%%%%%%%%%%%%%

\subsection{Nerve and Davis Complex}

Let $\Gamma$ be a right angled Coxeter group with generating set $S$.
The nerve $N=N(\Gamma, S)$ 
is the simplicial complex defined in the following way: the 
vertices of $N$ are the elements
of $S$. Two different vertices $s,t$ are joined by an edge, 
if and only if
$m(s,t)= 2$. In general (k+1) different vertices $s_1,\ldots ,s_{k+1}$ 
span a k-simplex, if and only
$m(s_i,s_k)=2$ for all pairs of different $i,j \in \{ 1,\ldots,,k+1\}$.
For a simplex $\sigma$ of $N$, let $\Gamma_{\sigma}$ be the 
subgroup of $\Gamma$ generated by the vertices of $\sigma$.
If $\sigma$ is a $k$-simplex spanned by
$s_1,\ldots,s_{k+1} \in S$ then $\Gamma_{\sigma}$ is
isomorphic to
$\Z_2^{k+1}$.
By $N'$ we denote 
the barycentric subdivision of $N$.
The cone $C=\cone N'$  over $N'$ is 
called a {\it chamber}
for $\Gamma$.
The Davis complex \cite{D} $X=X(\Gamma,S)$ is  
the image of
a simplicial map $q:\Gamma\times C \rightarrow X$ 
defined by 
the following equivalence relation on the vertices:
$a\times 
v_{\sigma}\sim b\times v_{\sigma}$ provided $a^{-1}b\in\Gamma_{\sigma}$.
Here 
$\sigma$ is a simplex in $N$, $\Gamma_{\sigma}$ is the subgroup of
$\Gamma$ 
generated by the vertices of $\sigma$, $v_{\sigma}$ is 
the 
barycenter of $\sigma$.
We identify $C$ with the image
$q(1\times C)$ as a subset of $X$.
The 
group $\Gamma$ acts simplicially on $X$ 
by $\gamma q(\alpha\times x) = q(\gamma\alpha\times x)$
and the orbit space is equal to the chamber $C$.
Thus the Davis complex $X$ is obtained
by 
gluing the chambers $\gamma C$ , $\gamma \in \Gamma$ 
along the boundaries.
Note that $X$ admits an equivariant cell structure with the 
vertices $X^{(0)}$
equal the cone points of the chambers and with 
the 1-skeleton $X^{(1)}$
isomorphic to the Cayley graph of 
$\Gamma$.

An alternative description of the Davis complex is obtained in
the following way:
Consider the cubical cell complex $\wh X$, whose
vertex set consists of the elements of $\Gamma$.
The 1-cells are of the form  $[\gamma,\gamma s]$ where $\gamma \in \Gamma$ and
$s\in S$. Thus the 1-skeleton $\wh X ^{(1)}$ is the Cayley graph $C(\Gamma,S)$.
The 2-cells are squares with a vertex set of the form
$\gamma,\gamma s,\gamma t, \gamma st =\gamma ts$, where 
$s,t \in S$ are distinct commuting
elements. In general the k-cells are cubes with vertex set
$\gamma \Gamma_{\sigma}$, where $\sigma$ is a $(k-1)$-simplex in $N$.
Then $\wh X$ can be considered as a cubical
realization of the Davis complex.

The generators $s\in S$ and their conjugates $r=\gamma s \gamma^{-1}$, 
$\gamma\in\Gamma$ are 
called {\it reflections}.
A {\it mirror (or wall)} of a 
reflection $r \in\Gamma$ is the set of 
fixed points $M_r\subset X$ 
of $r$ acting on the Davis complex $X$.
Note that $M_{\gamma s \gamma ^{-1}}= \gamma M_s$.

\begin{Lem} \label{mirror}
For every generator $s$ in a 
right-angled Coxeter group $\Gamma$ there is
the equality 
$M_s=\{ q(w \times x) \mid w\in Z_s(\Gamma), x\in St(s,N')\}$.
\end{Lem}

\begin{proof} 
$" \ \supset ": \ $ 
If $w\in Z_s(\Gamma)$ and $x\in St(s,N')$, i.e.
$x$ is an affine combination
$x = \sum_{s\in\sigma} x_{\sigma} v_{\sigma}$ one
easily computes
 $$s q(w\times x)=w q(s\times x)= w q(1\times x) = q(w\times x)$$.

\noindent $" \ \subset ": \ $ 
Let $z\in M_s$. 
Then $z=q(g\times x)$ for some $g\in\Gamma$ and $x\in cone(N')$. 
The 
condition $s(z)=z$  can be rewritten as $q(sg\times x)=q(g\times x)$.
Hence 
$g^{-1}sg\in\Gamma_{\sigma}$ 
for some simplex $\sigma$ of $N$ and $x\in\cap_{v\in\sigma} St(v,N')$. 
By 
the deletion law 
$s\in\Gamma_{\sigma}$, since the number of $s$ 
in $g^{-1}sg$ is odd.
Hence $s\in\sigma$ and $x\in St(s, N')$.

Let $s_1\dots s_k$ be a reduced presentation of $g^{-1}sg$.
We note that all $s_i\in\sigma$. Since the group is rightangled, $\Gamma_{\sigma}$
is commutative and hence all $s_i$ are different. Let $u_1\dots u_l$ be 
a reduced presentation of $g$. Note that $s_j\ne u_i$ for every $u_i\ne s$,
since $u_i$ appears even number times in the word $u_l\dots u_1su_1\dots u_l$.
Hence $g^{-1}sg=s$, i.e. $g\in Z_s(\Gamma)$.
\end{proof}

The {\it chromatic number} of a 
graph is the minimal number of colours
needed to colour the vertices 
in such a way that every adjacent vertices have 
different colours. A chromatic number of a simplicial complex
is 
the chromatic number of its 1-dimensional skeleton.

Assume that 
the chromatic number of the nerve $N(\Gamma)$ of a right-angled
Coxeter 
group $\Gamma$ equals $n$ and let $c:N^{(0)}\to\{1,\dots, n\}$
be a 
corresponding colouring map. Then for every mirror $gM_s$ 
in $X$ we can 
assign a colour by taking the colour $c(s)$.
Similarly we colour every edge
$[\gamma,\gamma s]$ of the Cayley graph
$C(\Gamma,S)$ with the colour $c(s)$.

\begin{Lem} \label{samecolour}

Different 
mirrors of the same colour are disjoint.
\end{Lem}

\begin{proof}
Let 
$c(s)=c(t)$, $s,t\in S$ and let $\gamma_1(M_s)\cap\gamma_2(M_t)
\ne\emptyset$.
Therefore $gM_s\cap M_t\ne\emptyset$ where $g=
\gamma_2^{-1}\gamma_1$. 
Let $x\in gM_s\cap M_t$. 
By Lemma \ref{mirror} we have
$x=q(w\times y)=q(gu\times z)$, where
$w\in Z_t(\Gamma)$, $y\in St(t, N')$
and
$u\in Z_s(\Gamma)$,$\ z\in St(s, N')$.
Hence $y=z = \sum_{t,s\in\Sigma} y_{\sigma} v_{\sigma}$.
In particular $s,t$ are in a common simplex and hence
commute. Since $c(s)=c(t)$, we see that $s=t$.
Since $q(w\times y)=q(gu\times z)$, we
have
$w^{-1}gu \in \Gamma_{\sigma}$ for some
simplex $\sigma$ with $s\in \sigma$.
Thus $w^{-1}gu \in Z_s(\Gamma)$ which implies
$g \in Z_s(\Gamma)$ and $gM_s=M_s$.

\end{proof}

Lemma \ref{samecolour} corresponds to the following
fact of the decomposition $\Gamma = H_s \cup s H_s$.

\begin{Lem}

$c(s) = c(t) \Rightarrow \gamma \partial H_s \subset H_t$ or
$\gamma \partial H_s \subset tH_t$

\end{Lem}

\begin{proof}
Assume first that $\gamma \in H_t$.
Let $a\in Z^0_s(S)$. If $\gamma a \in tH_t$ then
by Lemma \ref{halfspace} (b)
$\ell (\gamma a )> \ell (\gamma )$.
Let $s_1\ldots s_k$ be a reduced representation of 
$\gamma$, then $\gamma a$ is reduced.
Would $t s_1\ldots s_k a$ be not reduced
then by Lemma 
$t=a$, i.e. $t\in Z_a^0(S)$ and $c(t) \neq c(s)$
a contradiction.
Thus $t s_1\ldots s_k a$ is reduced and hence
$\ell ( t \gamma a )> \ell (\gamma a)$ and hence
$\gamma a \in H_t$.
By induction $\gamma \alpha \in H_t$ for
$\alpha \in Z_s^0(\Gamma)$ hence
$\gamma (\partial H_s)\subset H_t$.
If $\gamma \in tH_t$ then $t\gamma \in H_t$ and
by the above $t \gamma\partial H_s \subset H_t$
which implies $\gamma \partial H_s \subset t H_t$.

\end{proof}

%%%%%%%%%%%%%%%%%%%%%%%%%%%%%%%%%%%%%%%%%%%%%%%%%%%%%%%%%%%%%%%%%%%%%%%%%%%%%%%%%%%%%%%

\section{Maps into Trees}

In this section we study maps of
the Cayley graph into products of trees.

\subsection{Components of the Davis Complex} \label{components}

The mirror
$ M_s$ devides the Davis complex into two connected
components corresponding to the decomposition
$\Gamma = H_s \cup sH_s$.
We have 
$$X= (\bigcup_{\gamma\in H_s} \gamma C) 
\cup (\bigcup_{\gamma\in s H_s} \gamma C)$$
and the common boundary of this parts is $M_s$.
The two parts (and also the common wall ) are connected 
by Lemma \ref{halfspace} (d).

If we consider the Cayley graph $C(\Gamma,S)$ as the 1-skeleton
$X^{(1)}$ of the Davis complex, then every edge $e$ of the
Cayley graph intersects exactly one wall,
the edge $[\gamma,\gamma s]$ intersects the wall
$\gamma M_s$. Thus an edge $e$ with
colour $c(e)$ of the Cayley graph intersects a mirror with
the same colour.
We say that two edges
$e,e'$ of the Cayley graph are {\it parallel},
if $e$ and $e'$ intersect the same mirror.

\begin{proof} {\it (of Theorem \ref{tree})}

Let $c$ be a colour. Consider the graph
$T_c$ with vertices the connected components of 
$$X\setminus\bigcup_{c(s)=c,\gamma\in\Gamma}\gamma M_s$$
and edges correspond to the walls between components.
Since the Davis complex is simply connected \cite{D} and every wall devides $X$ in
exactly two pieces and different walls of the same color do not intersect, 
$T_c$ is indeed a tree (see \cite{DJ}).
We define a map $p_c:X^{(0)}\to T_c$ by the rule: $p_c(v)$ is the component
that 
contains $v$. This map extends simplicially to the Cayley graph 
$X^{(1)}$.

Note that $p_c$ is equivariant. 
The maps $p_c$ define an equivariant map
$\mu: C(\Gamma, S)\to\prod_{c=1}^nT_c$. 
Remember that we take the $l_1$-metric on
$\prod_{c=1}^nT_c$. Since every edge in the Cayley
graph intersects exactly one mirror, the distance between 
$\gamma_1$ and $\gamma_2$
in $\Gamma$ equals the number of walls between 
points
$q(\gamma_1\times x_0)$ and $q(\gamma_2\times x_0)$ where $x_0$ is 
the cone point
in the chamber $C$. On the other hand this number is exactly 
the 
$l_1$-distance in the product of our trees. Thus, $\mu$ is an 
isometry.
\end{proof}

%%%%%%%%%%%%%%%%%%%%%%%%%%%%%%%%%%%%%%%%%%%%%%%%%%%%%%%%%%%%%%%%%%%%%%%%%%%%%%%%%%%

\subsection{Locally Compact Trees}

In section \ref{components} we constucted a tree from the connected
components of the Davis complex. The edges of the $c$-tree correspond to
the $c$-mirrors or equivalently to the parallelclasses of $c$-edges of
the Cayley graph.
In this section we construct maps into locally compact trees.
In a certain sense the following construction is a modification of
the result in section \ref{components}. We have however to reformulate the result
in a different language:

We first introduce a general class of rooted simplicial trees:
Let $(Q)= Q_1,Q_2,\ldots $ be a sequence of nonempty sets.
We associate to $(Q)$ a rooted simplicial tree $T_{(Q)}$ in the following way
(compare \cite{dH} p. 211):
The set of vertices is the set of finite sequences
$(q_1,\ldots ,q_k)$ with
$q_{\kappa} \in Q_{\kappa}$.
The empty sequence defines the root vertex and is denoted
by $v_{\emptyset}$.
We write the vertex given by $(q_1,\ldots ,q_k)$ also
as $v_{(q_1,\ldots ,q_k)}$.
Two vertices are connected by an edge in
$T_{(Q)}$ if their length (as sequences) differ by one
and the shorter can be obtained by erasing the last term of the longer.
The root vertex has $\mid Q_1 \mid$ neighbours and every vertex of distance
$i$ to $v_{\emptyset}$ has
$\mid Q_{i+1} + 1\mid$ neighbours, one ancestor and
$\mid Q_{i+1} \mid $ descendents. Here $\mid Q_{i} \mid$ denotes the
cardinality of $Q_i$.
The tree $T_{(Q)}$ is locally compact, if and only if
$Q_i$ is finite for all $i$.

We recollect and extend certain notations.
We have given a colouring
$c: S \to \{1,\ldots, n\}$ such that
$c(s) \neq c(t)$ if $ts = st$ and $t \neq s$.
We have the length function
$\ell : \Gamma \to \N$ and we
introduce for a colour $c$ also the
function
$\ell_c : \Gamma \to \N$ by defining
$\ell_c(\gamma)$ to be the number of letters with
colour $c$ in a reduced representation of $\gamma$.
This does not depend on the special representation.
We denote by
$\cal R  \subset$ $ \Gamma $
the set of all reflections.
An element $r \in \cal R $ can be represented as
$r = \gamma s \gamma^{-1}$ for some $s \in S$ and some
$\gamma \in \Gamma$.
One easily checks that
$\gamma_1 s \gamma_1^{-1} = \gamma_2 t \gamma_2^{-1}$ for
$\gamma_i \in \Gamma$ and
$s,t \in S$ if and only if $s=t$ and $\gamma_1^{-1}\gamma_2 \in Z_s(\Gamma)$.
This implies that the maps
$g:{\cal R}  \to S$, $\gamma s \gamma^{-1} \mapsto s$ 
and
$c:{\cal R}  \to \{ 1,\ldots , n\}$, $c(r) = c(g(r))$ are
well defined.

As defined above the mirror $M_r$ of a reflection $r \in \cal R$ 
is a subset of the Davis complex. It intersects the Cayley graph
$C(\Gamma , S)$ (viewed as embedded in the Davis
complex) in a set of parallel edges:
A given edge
$ e = [\gamma,\gamma s ]$ of the Cayley graph
intersects the mirror
$M_r$ of the reflection
$r(e) = \gamma s \gamma^{-1} \in \cal R$.
If on the other side a reflection $r \in \cal R$ can be represented as
$ r = \gamma s \gamma^{-1}$ then $r$ interchanges $\gamma$ and
$\gamma s$ and hence $M_r$ intersects the edge of $C(\Gamma ,S)$ with
endpoints $\gamma $ and $\gamma s$. By replacing $\gamma$ by $\gamma s$ if 
necessary we can assume that
$\ell (\gamma s) > \ell (\gamma)$ and hence we can write this
edge as the oriented edge
$[\gamma , \gamma s]$ .
We define the level function 
$\lev : \cal R \to \N$ by
$\lev (r) = \ell_c(\gamma) + 1$, where 
$ c = c(r)$ and $ \gamma$ is an element such
that $r = \gamma s \gamma^{-1}$ and $\ell (\gamma s) > \ell (\gamma)$.
It is not difficult to show that the level function is well defined.
The level function can also be considered as a function defined on the
set of edges of the Cayley graph.
This function has an easy geometric interpretation:
Let $ r \in \cal R$ with $c = c(r)$. Consider a shortest path in the 
Cayley graph from the origin $1$ to the mirror $M_r$. Then
$\lev (r)$ is the number of mirrors with colour $c$ which are 
intersected by this path
(including the final mirror $M_r$).

We denote by
${\cal R}^c \subset \cal R$ the reflections with colour $c$,
and by ${\cal R}^c_i$ the set of reflections with colour $c$ and
level $i$. We consider the tree
$T_{({\cal R}^c)}$ belonging to the sequence
$ ( {\cal R}^c ) = {\cal R}^c_1,{\cal R}^c_2,\ldots $.

For a colour $c$ we define
a map $\phi_c:\Gamma \to T_{({\cal R}^c)}$ in the following
way:
Let $\gamma \in \Gamma$ and let $s_1\cdots  s_k$ be a reduced word representing
$\gamma$.
This corresponds to a geodesic path following the edges
$e_1,\ldots ,e_k$ in the Cayley graph, where
$e_i = [s_1\cdots s_i, s_1\cdots s_{i+1}]$.
Let $1\leq i_1 < \ldots < i_{\ell_c(\gamma)} \leq k$ be the set of
indices with
$c(s_{i_{\kappa}})=c(e_{i_{\kappa}}) = c$.
Define
$$ r_{\kappa} =
r(e_{i_{\kappa}}) = s_1\cdots s_{i_{\kappa -1}} s_{i_{\kappa}} 
(s_1\cdots s_{i_{\kappa -1}})^{-1} \in {\cal R}^c $$
and
$$\phi_c(\gamma) = (r_1,\ldots ,r_{\ell_c(\gamma)})$$
By construction
$\lev (r_{\kappa}) = \kappa$ and $\phi_c(\gamma) \in T_{({\cal R}^{c})}$.

We show that $\phi_c$ is well defined.
By the deletion rule we have to check the following:
let $w_1stw_2$ and $w_1tsw_2$ are reduced representations
of the same element $\gamma$, then the above definition gives the same
$\phi_c$-image for both words.
Note that $s \neq t$ and $s$ commutes with $t$. Thus $s$ and $t$
have different colours.
If $s$ and $t$ both do not have the colour $c$ then the sequence
$(r_1,\ldots ,r_{\ell_c(\gamma)})$ does not change at all.
If one of $s$ or $t$ has colour $c$ we assume w.l.o.g. that
$c(t) = c$ and $c(s) \neq c$.
Let this $t$ be the $\kappa$-th letter of colour $c$ in the word
$w_1stw_2$, then it is also the $\kappa$-th letter of colour $c$ in the word
$w_1tsw_2$.
Clearly the reflections 
$r_1,\ldots,r_{\kappa -1},r_{\kappa +1},\ldots,r_{\ell_c(\gamma)}$
are the same for the two words.
Since
$ r_{\kappa} = w_1st(w_1 s)^{-1} = w_1 t w_1^{-1}$ 
we see that $r_{\kappa}$ is also the same
in both words we are done.

\begin{Rem}
The tree $T_c$ from section \ref{components} is canonically isometrically
embedded into $T_{({\cal R}^c)}$ by the following map $\chi$:
The vertices of $T_c$ are the components of the Davis complex minus the $c$-walls.
The component containing the base Chamber $C$ is mapped by $\chi$ to
the root vertex $v_{\emptyset} \in T_{({\cal R}^c)}$. If
$C'$ is some other component, consider a shortest path from
$C$ to $C'$. This path intersects a sequence of $c$-walls
$M_{r_1},\ldots,M_{r_k}$. Now define
$\chi (C') = v_{(r_1,\ldots,r_k)}$.
Clearly the image of $\phi_c$ is contained
in $\chi (T_c)$. By identifying $T_c$ with its image we
obtain a map
$\phi_c: \Gamma \to T_c \subset T_{({\cal R}^c)}$. This map
is exactly the map $p_c$ defined in section \ref{components}.

\end{Rem}

The sets ${\cal R}^c_i$ are (besides of trivial cases) not finite and hence the
corresponding tree is not locally compact.
In order to get a map into locally compact trees we
have to replace the infinite sets ${\cal R}^c_i$ by certain finite
sets. Therefore we will construct
maps ${\fin}^c_i:{\cal R}^c_i \to F^c_i$ where $F^c_i$ are finite sets.
These maps have to satisfy the condition that two different 
reflections whose mirrors are close to each other
are mapped to different points.
The distance between mirrors is defined as follows

\begin{Def}

Let $M_{r_1}$ and $M_{r_2}$ be two different mirrors with the same colour
$c=c(r_1)=c(r_2)$.
Then the distance is defined to be the number
$$d(M_{r_1},M_{r_2}) = \inf \{ d(\gamma_1,\gamma_2) + 1 \mid
\gamma_i \in \Gamma,,\ s_i \in S, \ r_i = \gamma_i s_i \gamma_i^{-1}\} $$

\end{Def}

\begin{Rem}

We can view the mirrors $M_{r_i}$ as subsets of the Cayley graph such
that $M_{r_i}$ is the set of midpoints of the edges with
endpoints $\gamma_i $ and $\gamma_i s_i$ where
$\gamma_i s_i \gamma_i^{-1}$ is a representation of $r_i$.
Since by Lemma \ref{samecolour} 
two different mirrors with the same colour do not
intersect, the distance defined above is exactly the distance
of the mirrors considered as subsets of the Cayley graph.

\end{Rem}

The following result is essential for our construction

\begin{Pro} \label{Proposition}

There exists a map
${\fin}^c_i:{\cal R}^c_i \to F^c_i$ where $F^c_i$ is a finite
set such that
${\fin}^c_i (r_1) = {\fin}^c_i (r_2)$ implies
$r_1 = r_2$ or $d(M_{r_1},M_{r_2}) \geq 4 n i$.
Furthermore there exists a constant $\rho > 0$ (independent
of $i$) such that
$\mid F^c_i \mid \leq e^{\rho i}$.

\end{Pro}

In order to prove Proposition \ref{Proposition} we need

\begin{Lem} \label{finitegroup}

Let $\nu \in \N$ be given.
Then there exists a finite group $F_{\nu}$ and a
homomorphisms
$\sigma : \Gamma \to F_{\nu}$
with the property:
If $r_1,r_2 \in {\cal R}$ with 
$g(r_1) = g(r_2)$,
$d(M_{r_1},M_{r_2}) < \nu$ and
$\sigma (r_1 ) = \sigma (r_2)$
then
$r_1 =r_2$.
Furthermore there exists a constant
$\rho_1$ (independent of $\nu$) such that
$\mid F_{\nu}\mid \leq e^{\rho_1\nu}$.

\end{Lem}

\begin{proof}
Let 
$\alpha_1,\ldots ,\alpha_k \in \Gamma$ be the
set of nontrivial elements 
with
$\ell (\alpha_i) \leq 2\nu +2 $.
Since a Coxeter group $\Gamma$ is
residually finite, there exists a finite
group
$F_{\nu}$ and a homomorphism
$\sigma : \Gamma \to F_{\nu}$
such that
$\sigma (\alpha_i) \neq 1$ for all
$i=1,\ldots,k$.
Let $r_1,r_2 \in \cal R$ with
$g(r_1) = g(r_2)$ and
$d(M_{r_1},M_{r_2}) < \nu$.
Then there exits $s \in S$ and $\gamma_j \in \Gamma$ such that
$r_j = \gamma_j s \gamma_j^{-1}$ and
$d(\gamma_1,\gamma_2) \leq \nu$.
Let
$\tau = \gamma_1^{-1}\gamma_2$,
hence $\ell(\tau) \leq \nu$.
By assumption
$r_1^{-1}r_2 \in ker(\sigma)$.      
Since
$r_1^{-1}r_2=\gamma_1s\tau s \tau^{-1}\gamma_1^{-1}$
we have
$s \tau s \tau^{-1} \in ker (\sigma)$.
Since
$\ell(s \tau s \tau^{-1}) \leq 2 \nu + 2$
we have by construction that
$s \tau s \tau^{-1}$ is trivial and
hence
$\tau$ commutes with $s$ which implies that
$r_1 = r_2$.

We finally have to estimate the size of $F_{\nu}$.
Therefore we use a concrete geometric realisation of
$\Gamma$ as a group of linear transformations on a vectorspace
$V$ over $\R$ (compare \cite{H} p. 108).
$V$ has the basis $\{ v_s \mid s\in S\}$ in one to one
correspondence to $S$.
For each $s \in S$ define a reflection
$h(s): V \to V$ by
$h(s)(v_s) = - v_s$,
$h(s)(v_t) = v_t + 2 v_s$ if $t \neq s$ and $ts \neq st$,
$h(s)(v_t) = v_t - 2 v_s$ if $t \neq s$ and $ts = st$.
Then $h$ defines a faithful representation of $\Gamma$ and
enables us to identify $\Gamma$ with a subgroup of
$GL(\mid S\mid,\Z)$.
The description also implies that every matrix coefficient of
$h(\alpha)$ is an integer with norm
$\leq 3^{\ell(\alpha)}$.
Let $F_{\nu} = GL(\mid S\mid, \Z/(3^{2\nu +2}+1) \Z)$,
then the canonical map
$\sigma : \Gamma \to F_{\nu}$ satisfies
$\sigma (\alpha ) \neq 1$ if
$0 < \ell (\alpha ) \leq 2\nu +2$.
Clearly $\mid F_{\nu} \mid\  \leq e^{\rho_1 \nu}$ where 
$\rho_1$ only depends on $\mid S \mid$.

\end{proof}

\begin{proof} (of Proposition \ref{Proposition})

Let $c \in \{ 1,\ldots , n\}$ be a given colour and $i \in \{1,2,\ldots\}$ be
a given level.
Let $S^c \subset S$ be the subset of generators with colour $c$.
For $s \in S^c$ let
${\cal R}^s_i$ be the set of reflections with generator $s$ and level $i$,
${\cal R}^c_i = \bigcup_{s\in S^c}{\cal R}^s_i$.
By Lemma \ref{finitegroup} there exists a set
$F^s_i$ and a map
$\fin^s_i:{\cal R}^s_i \to F^s_i$ with the property:
if $\fin^s_i(r_1) = \fin^s_i(r_2)$ then
$r_1 = r_2$ or
$d(M_{r_1},M_{r_2}) \geq 4 n i$.
Define $\fin^c_i$ to be the sum of the maps $\fin^s_i$ , $s \in S^c$ from
${\cal R}^c_i$ to $F^c_i$.

\end{proof}

\begin{Rem}

For simplicity we use the notation
$\fin$ instead of $\fin^c_i$ if the indices are clear from the context.

\end{Rem}

For a colour $c$ we consider the locally compact tree
$T_{(F^c)}$ coming from the sequence
$(F^c)=F^c_1,F^c_2,\cdots $.
We consider the map
$\psi_c : \Gamma \rightarrow T_{(F^c)}$
defined by
$$\psi_c(\gamma) = (\fin (r_1),\ldots, \fin (r_{\ell_c{(\gamma)}}))$$
where
$$\phi_c(\gamma) = (r_1,\ldots ,r_{\ell_c(\gamma)})$$

This map extends naturally to a map
$$\psi_c:C(\Gamma,S)\to T_{(F^c)}$$
We define
$$\psi = \prod_{c=1}^n \psi_c: 
C(\Gamma,S)\rightarrow \prod_{c=1}^n T_{(F^c)}$$

Essential for the bilipschitz property of the map
$\psi$ is the
following result

\begin{Lem} \label{mainlemma}
Let $\gamma , \gamma' \in \Gamma$
and let
$d(\psi(\gamma),\psi(\gamma'))=m$
then
$d(\gamma,\gamma')\leq 16 n m$.
\end{Lem}

\begin{proof}

By Lemma \ref{between}
there are 
exists an element
$\alpha \in \Gamma$ between
$\gamma$ and $\gamma'$ such that
$\alpha$ lies also between $1$ and $\gamma$
and $1$ and $\gamma'$.
We now consider a geodesic
$\alpha=\alpha_0,\cdots,\alpha_{\tau}=\gamma$
from
$\alpha $ to $\gamma$ and a geodesic
$\alpha =\alpha_0',\cdots,\alpha'_{\tau'}=\gamma'$ from
$\alpha$ to
$\gamma'$.
Then $\tau+\tau' = d(\gamma,\gamma')$. 
The edges
$e_i=[\alpha_{i-1},\alpha_i]$ and
$e_{i'}=[\alpha_{i-1}',\alpha_i']$ are then
oriented edges of the Cayley graph and the path
$e_{\tau},\ldots,e_1,e_1',\ldots,e_{\tau'}'$ is a geodesic in 
$C(\Gamma,S)$ from $\gamma$ to $\gamma'$ with length
$\tau+\tau'= d(\gamma,\gamma')$.  We can assume without loss of
generality that
$\tau \geq \tau'$.
Let $\tau_c$ be the number of $c$-edges in the path $e_1,\ldots,e_r$.
Choose the colour $c\in \{ 1,\ldots,n\}$ in a way that
$\tau_c$ is maximal.
If $\tau_c \leq 8m$ then
$d(\gamma,\gamma') = \tau + \tau' \leq 2n \tau_c \leq 16 n m$
and we are done.

Thus we can assume $\tau_c > 8m$. Let $e_j$ be the 
$(\tau_c-m)$-th $c$-edge in the geodesic path
$e_1,\ldots, e_{\tau}$ and let
$r_j$ be the corresponding reflection.
Since  
$m = d(\psi(\gamma),\psi(\gamma'))\geq d(\psi_c(\gamma),\psi_c(\gamma'))$
there exists an $c$-edge
$e'$ in the path $e_1',\ldots ,e_{\tau'}'$ with 
corresponding reflection $r'$
such that
$\lev(r')=\lev(r_j)$ and $\fin (r')=\fin (r_j)$.
We claim that $r_j = r'$ .
Note that $\lev(r_j) \geq (\tau_c-m)$ and
$d(M_{r_j},M_{r'}) < \tau+\tau' \leq 2 n \tau_c$.
If $r_j \neq r'$, then by Proposition \ref{Proposition}
we have

$$d(M_{r_j},M_{r'})\geq 4 n \lev(r_j)) \geq 4 n (\tau_c-m) 
\geq 4 n (\tau_c -\frac{\tau_c}{8}) = 2 n \tau_c$$
a contradiction.

Thus $r_j = r'= \beta s \beta^{-1}$
for some $\beta \in \Gamma$ and some $s\in S$.
The edges $e_j= [\alpha_{j-1},\alpha_{j-1} s]$ and 
$e' = [\alpha_{j'-1}',\alpha_{j'-1}'s]$ are parallel and
intersect the same mirror
$\beta M_s$.
Thus $\alpha_{j-1}$ and $\alpha_{j'-1}'$ are both in $\beta H_s$ and
hence by Lemma \ref{halfspace}
the geodesic
$\alpha_{j-1},\ldots ,\alpha_0=\alpha_0',\ldots,\alpha_{j'-1}'$ is completely
contained in $\beta H_s$ and in particular
$\alpha \in \beta H_s$.
Now $\alpha_{j}$ and $\alpha_{j'}'$ are contained
in $\beta sH_s$ and by the same argument the complete geodesic
$\alpha_{j},\ldots ,\alpha_0=\alpha_0',\ldots ,\alpha_{j'}'$
is contained in $\beta sH_s$. Hence 
$\alpha \in \beta H_s \cap \beta sH_s = \emptyset$.
Thus the assuption $\tau_c > 8 m$ leads to a contradiction.
\end{proof}

\begin{proof} (of Theorem \ref{locallycompacttree})

The map
$\psi:C(\Gamma,S) \to T_{(F^c)}$ is clearly
1-lipschitz. By Lemma \ref{mainlemma} $\psi$ is bilipschitz.
Let $T$ be any exponentially branching
tree. Since the tree $T_{(F^c)}$ satisfies the estimate
$\mid F^c_i\mid \leq e^{\rho i}$, one can show that there
exists a bilipschitz embedding of $T_{(F^c)}$ into $T$
(we do {\it not} require that the root vertex of $T_{(F^c)}$ is
mapped to a given basevertex $t_0 \in T$).
Combining these results we obtain a bilipschitz embedding of
$C(\Gamma,S)$ into $T_{(F^c)}$.
Since the maps
$\fin^c_i:({\cal R}^c_i,d) \to F^s_i$ 
where 
$d(R_1,r_2) = d(M_{r_1},M_{r_2})$ are locally
injective by \ref{Proposition}, the map
$\psi$ has locally the same properties as the map
$\mu$, i.e. it is locally an isometry.
By adjusting the constants in \ref{Proposition} suitable,
we can enforce that $\psi$ is an isometric embedding on every ball
of a given radius $r$.
\end{proof}

\begin{proof} (of Corollary \ref{hypplane})

Consider the right angled Coxeter group
$\Gamma$ given by the generator set $S=\{s_1,\ldots,s_6\}$
and relations $s_i s_{i+1} = s_{i+1} s_i $ (indices mod 6).
This group acts discretely on the hyperbolic plane 
such that a Dirichlet fundamental
domain is bounded by the regular right angled hexagon in $\H^2$.
By
Theorem \ref{locallycompacttree}.
we can embed the Cayley graph of $\Gamma$ locally isometric
and globally bilipschitz into a product of trees.
The Cayley graph can be realized canonically in $\H^2$, where
the vertices are the central points of the hexagonal decomposition
of $\H^2$ and the edges are geodesics. Then the Cayley graph gives
a decomposition of $\H^2$ into regular quadrilaterals (with angles
equal to $\frac{\pi}{3}$). The map $\psi$ maps the boundary of
a square to a square in $T\times T$. Thus the maps $\psi$ can be extended
from the Cayley graph in a bilipschitz way to all of $\H^2$.
\end{proof}

\begin{proof} (of Corollary \ref{higherdim})

By a result of Brady and Farb \cite{BF} there exists a bilipschitz embedding of
the hyperbolic space $\H^n$ into the $(n-1)$-fold product of hyperbolic
planes. Actually in \cite{BF} it is only stated that the 
embedding is quasiisometric
but their proof gives a bilipschitz embedding 
(compare \cite[section 2]{F}).
Combining Corollary \ref{hypplane} with this result, we are done.

\end{proof}

\begin{proof}  (of Corollary \ref{vcohdim})

A recent result of Januszkiewicz and Swiatkowski \cite{JS}
shows the existence of a Gromov-hyperbolic right angled 
Coxeter group $\Gamma_n$
of arbitrary given n, such that the virtual cohomological 
dimension of $\Gamma_n$ is 
$n$. The construction of these groups imply that they have chromatic number
$n$. By Theorem \ref{locallycompacttree} $\Gamma_n$ can be embedded in a
bilipschitz way into the product $T^n$.

\end{proof}

\begin{proof} (of Corollary \ref{hyprank})

We recall the definition of the hyperbolic rank of a metric space.
Given a metric space $M$ consider all locally compact Gromov-hyperbolic
subspaces $Y$ quasiisometrically embedded into $M$. Then
$\rank_h(M) = \sup_Y \dim \partial_{\infty}Y$ is called
the hyperbolic rank. (Compare \cite{BS1} for a discussion of this
notion). 
Let $T$ be an exponentially branching tree.
Let $\Gamma_n$ as in the proof of Corollary \ref{vcohdim} above.
$\Gamma_n$ can be embedded in a
bilipschitz way into the product $T^n$. Since 
the virtual cohomological dimension of $\Gamma_n$ is $n$, we have
$\dim \partial_{\infty} \Gamma_n = (n-1)$ 
by  \cite{BM}.
Thus $\rank_h(T^n) \geq (n-1)$ by the definition of the hyperbolic rank.
The opposite inequality
$\rank_h(T^n) \leq (n-1)$ follows from standard topological considerations.
 
\end{proof}

%%%%%%%%%%%%%%%%%%%%%%%%%%%%%%%%%%%%%%%%%%%%%%%%%%%%%%%%%%%%%%%%%%%%%%
%%%%%%%%%%

\bigskip
\begin{tabbing}

Alexander Dranishnikov,\hskip10em\relax \= Viktor Schroeder,\\ 

Dep. of Mathematics,\>
Institut f\"ur Mathematik, \\

University of Florida,\> Universit\"at Z\"urich,\\
444 Little Hall, \>
 Winterthurer Strasse 190, \\

Gainesville, FL 32611-8105\>  CH-8057 Z\"urich, Switzerland\\

{\tt dranish@math.ufl.edu}\> {\tt vschroed@math.unizh.ch}\\

\end{tabbing}
\end{document}